\input amstex
\documentstyle{amsppt}
\magnification=\magstep1
\define\cc{\Bbb C}
\define\r{\Bbb R}

\define\Z{\Bbb Z}
\define\ap#1{\Cal A\Cal P(#1)}
\define\aap#1{\Cal A\Cal A\Cal P(#1)}
\define\as#1{\Cal A\Cal A(#1)}
\define\pas#1{\Cal P\Cal A\Cal A(#1)}
\define\C#1{\Cal C(#1)}
\define\co#1{\Cal C_0(#1)}
\define\cu#1{\Cal C_u(#1)}
\define\pap#1{\Cal P\Cal A\Cal P(#1)}

\define\wap#1{\Cal W\Cal A\Cal P(#1)}
\define\wapo#1{\Cal W\Cal A\Cal P_0(#1)}
\define\wpapo#1{\Cal W\Cal P\Cal A\Cal P_0(#1)}
\define\papo#1{\Cal P\Cal A\Cal P_0(#1)}

\define\ga{\gamma }
\define\va{\varphi }

\define\la{\lambda}
\define\La{\Lambda}

\define\st{\subset}
\define\alf{\alpha}
\topmatter
\title
Beurling Spectra        of                Functions         on
Locally Compact Abelian Groups
\endtitle
\author
B. Basit and A.J. Pryde
\endauthor
\abstract{Let  $G$   be a locally compact abelian  topological
group.   For  locally
bounded  measurable  functions $\varphi:  G\to\Bbb  {C}$   we  discuss notions  of
spectra  for  $\varphi$ relative  to  subalgebras  of
$L^{1}(G)$.  In particular we study polynomials on  $G$
and  determine  their  spectra.   We  also  characterize   the
primary   ideals   of certain Beurling algebras $L_{w}^{1}(\Bbb Z)$
on the group of integers  $\Bbb Z$.  This allows  us  to
classify those elements of $L_{w}^{1}(G)$  that
have finite spectrum.  If  $\varphi$  is  a  uniformly  continuous function  whose
differences are bounded, there is  a  Beurling algebra  naturally associated  with
$\varphi$.  We give a condition on the  spectrum of  $\varphi$  relative  to this  algebra
which ensures that  $\varphi$  is bounded.     Finally we   give  spectral conditions on a
bounded function on   $\Bbb R$ that   ensure   that   its   indefinite integral  is  bounded.}
\endabstract
\endtopmatter
\leftheadtext{B. Basit and A.J. Pryde}
\rightheadtext{Spectra of Unbounded Functions on Groups}
\document
\pageno=1
\baselineskip=20pt

\head{\bf\S1. Introduction, notation and preliminaries}\endhead
Throughout  this  paper  $G$  will denote  a  locally  compact
abelian   group   equipped with Haar measure.
We  consider certain unbounded functions   $\va:  G
\to\cc$   and  discuss
various  notions of spectrum for $\va$.  Our aim is to  deduce
properties  of   $\va$
from a knowledge of its spectrum and conversely.

Firstly, the spaces  $L^{p} (G)$  with corresponding norms $\|.\|_{p}$   are defined
as  usual  using  the  Haar  measure  (see [20]).    We   will
always  identify
functions   on   $G$   which  agree  almost  everywhere (see  [14, p. 36]).    In
particular,  $L^{1} (G)$   is a
Banach  algebra under the  convolution  product  and   carries
an  involution defined by

(1.1)  $f^{*} (x)=\overline {f(-x)}\text {    for    }  x\in G,  f\in  L^{1}(G).$

By    $L^{\infty}_{loc}(G)$     we   mean   the    space    of
functions   $\va:G\to\cc$   such   that
$\chi_{K}   \va\in  L^{\infty}(G)$   for  every  compact   set
$K\subset G$,  $\chi_{K}$   denoting the  characteristic
function of  $K$.

The action of the dual group  $\widehat {G}$  on  $G$  will be
denoted   $\ga (x)=(x,\ga)$
where  $  x\in  G$,  $\ga\in\widehat  G$   and   the   Fourier
transform   of  $f\in L^{1}(G)$  is defined  by
$\hat   {f}(\ga)=\int_{G}  f(x)(-x,\ga)\,dx$   for     $\ga\in
\widehat {G}$.   Unless  otherwise  stated,  the  group
operation  on   $G$   will be addition and on  $\widehat  G$
multiplication.  The unit in  $\widehat G$
is then  $\hat e$.
  The translate  $\va_{y}$   and  difference   $\Delta_{y}\va$
by  $y\in G$   of  a  function
$\va: G\to\cc$  are defined by

(1.2) $\va_{y}   (x)=\va(x+y)\text  {  and  }  \Delta_{y}\va(x)=\va
(x+y)-\va (x) \text { for } x\in G.$

 In  what  follows,   $F$   will  denote  a  dense  subalgebra
of $L^{1}(G)$  satisfying

(1.3) $F$  is translation invariant,

(1.4)     $F$    is  closed under pointwise  multiplication  by
characters  $ \ga\in\widehat {G},$

(1.5)     for each  $ \ga\in \widehat {G}$,  the ideal of  $ F$
   defined  by $ I(\ga,F) = \{f \in F:\hat  {f}(\ga) =0\} $   is dense  in
the  maximal ideal of $L^{1}(G)$ defined by
  $I(\ga) = \{f \in L^{1}(G):\hat {f}(\ga) = 0\}.$

      For  $f\in F$  and $\va\in L^{\infty}_{loc}(G)$  we will
say that the convolution   $f^{*} *\va$
is     well-defined   if    $f_{y}\va$  and    $f\va_{y}   \in
L^{1}(G)$    for   almost every    $y\in G$.    Then
$f^{*} *\va (x)=\int_{G} f^{*} (x-y)\va(y)\,dy=\int
f^{*}(y)\va (x-y)\,dy$  almost everywhere  (a.e)   in
$G$.

  If   $\va\in  L^{\infty}_{loc}(G)$  and   $f^{*}  *\va$   is
well-defined for all  $f\in F$  then
       $I_{F} (\va) = \{f \in F:  f^{*} *\va = 0\}$
is an ideal of the algebra  $F$.  We define the spectrum of
 $\va$  relative to  $F$ to be the set

(1.6) $   sp_{F} (\va)=\text {hull }I_{F} (\va)=\{ \ga \in \widehat {G}:\hat
{f}(\ga)=0 \text{ for all  }  f\in I_{F} (\va)\}.$

\noindent As  each such function  $\hat {f}$  is continuous on $\widehat
{G}$, $sp_{F} (\va)$  is a  closed  subset
of  $\widehat {G}$.

 The paper is written in five sections.  In section 2, we firstly
give the  basic
properties of the general spectrum  $sp_{F} (\va)$. Secondly,
  we   deal  with  a  large
class of examples, namely the Beurling  algebras  $F=L^{1}_{w}
(G)=L^{1}(\mu)$   where
$d\mu  (x)=w(x)dx$   is  a weighted  Haar   measure   on   $G$.
These  algebras  were
introduced  by Beurling [2] for the case  $G=\r$  and  studied
in  the  general
case  by  Domar [5]. Finally,
 we consider a smaller  dense subalgebra
of   $L^{1} (G)$,
namely   $F=C_{c} (G)$  the space of continuous complex valued
functions  on   $G$
with  compact  support.
      We  give  a new definition of polynomials on   $G$,   in
section 3, and prove
the  equivalence of our definition with one of Domar [5].   If
$\va =\ga p$   where
$\ga \in \widehat {G}$  and  $p$  is a polynomial on  $G$,  we
show  that $sp_{c} (\va)=\{\ga \}$.
     In section 4  we  again  consider  spectra  with  respect
to  Beurling
algebras,  but with additional conditions on the weights  $w$.
In  particular
we   describe  explicitly  the  closed  primary   ideals    of
$L^{1}_{w}(\Z)$    with
cospectrum  $\{ \hat {e}\}$.  This enables us to  characterize
those  $\va \in L^{\infty}_{w} (G)$  for which
$sp_{w}  (\va)$    is   finite.   They   are   precisely   the
generalized  trigonometric
polynomials    $\va=\sum^{n}_{k=1}  \ga_{k}   p_{k}$     where
$\ga_{k}\in \widehat G$   and $p_{k}$    is  a  polynomial  in
$L^{\infty}_{w} (G)$. In section 5 we consider  uniformly  continuous functions
$\va: G\to \cc$
whose differences  $\Delta_{y}\va, y\in G$,  are all  bounded.
Such  functions  have  a
corresponding weight function  $w_{\va}$   with  which  we may
define  a  spectrum
$sp_{w_{\va}}  (\va)$.

 Before ending this introduction, we make a few historical
comments  on
spectra.   For $\va\in L^{\infty} (\r)$,  Beurling  introduced
in [3] and [4] the spectral  set
$\La_{\va}  =\hat {\r}\cap \overline {[\va]}_{w^{*}}$,     where
$\overline  {[\va]}_{w^{*}}$  is the  weak-star   closure   of
the  smallest
translation   invariant   subspace    of   $L^{\infty}   (\r)$
containing $\va$.   Godement [9]
extended the definition of  $\La_{\va}$   from  $\va\in
 L^{\infty}(\r)$  to $\va\in L^{\infty}(G)$   and  proved
that   $\La_{\va}  =sp_{F}(\va)$, where $F=L_{1}(G)$.  The approach  of  Godement  was
pursued by Kaplansky [11] to
find   the   elements  of   $L^{\infty} (G)$   whose  spectrum
is  a singleton,  and  by
Helson[ 10] to give a new proof of Kaplansky's result.
      Domar  [5, ch.1]  introduced  an  abstract  setting  for   the
definition of spectra,
using  a  class of Banach algebras  $F$  more   general   than
Beurling  algebras.
However,  for  a  given   $\va\in  L^{\infty}_{loc}(G)$,   the
problem of finding an algebra $F$  for
which  $sp_{F}(\va)$  is defined was not discussed.

\head{\bf\S2. Spectra relative to subalgebras}\endhead

      In  this section we establish the main properties of the
spectra  $sp_{F} (\va)$
relative  to  subalgebras   $F$  of   $L^{1}(G)$ satisfying (1.3)-(1.5).
  We also give examples of such subalgebras.    For   $F=L^{1}(G),$ most
parts  of the following theorem can be found scattered in  the
literature ([5], [6], [19], [20]).  We include them here for reference purposes.

  \proclaim{Theorem 2.1} Let  $F$  and $F_{1}$  be    dense  subalgebras of
$L^{1}(G)$   satisfying
(1.3)-(1.5).   Let   $\va,\psi \in L^{\infty}_{loc}(G)$    and
suppose  $f^{*} *\va$ and  $f^{*} *\psi$    are
well-defined for each  $f\in F$.  Then

(a)  $sp_{F} (\va_{y} ) = sp_{F} (\va)$  for all  $y \in
 G$.

(b)   $sp_{F}  (k\va) = sp_{F} (\va)$  for  all   $k  \in\cc,
k\not=0$ .

(c)  $sp_{F} (\va+\psi)\subset sp_{F} (\va)\cup sp_{F}
(\psi)$.

(d) $sp_{F} (\ga \va)= \ga sp_{F} (\va)$  for all $\ga \in
\widehat G$.

(e)  $sp_{F}(\Delta_{y}\va) \subset sp_{F}  (\va)\subset \{\hat{e}\}\cup_{z\in G}
sp_{F}(\Delta_{z}\va)$  for all  $y \in G$.

(f)  $sp_{F}  (\va)   =\{\hat {e}\}$  if   $\va$   is  a  non-zero
constant, and $sp_{F} (\va)  = \emptyset$ if $\va=0$.

(g)  $sp_{F} (\va)\subset sp_{F_{1}} (\va)$  if  $F_{1}\subset
F$.
\endproclaim
\demo{Proof}  Since   $I_{F}(\va_{y})=I_{F} (k\va)=I_{F}(\va)$
for $y\in G$   and $k\not=0$,   (a)
and  (b)  follow.  For (c), assume  $\ga\not\in  sp_{F}  (\va)
\cup sp_{F} (\psi)$.   Then  there  exist
$f\in   I_{F}   (\va) $ and   $g\in   I_{F}   (\psi)$    such   that
$\hat{f}(\ga)\not=0$   and $\hat{g}(\ga)\not=0$. If $h=f*g$
then     $h\in   I_{F}   (\va+\psi)$    and,   since     $\hat
{h}(\ga)\not=0$,  we conclude $\ga\not\in sp_{F}  (\va+\psi)$.
Next,  if   $f\in F$  and  $\ga\in \widehat G$   then   $f^{*}
*\ga \va=\ga (\ga^{-1} f^{*}*\va)=\ga ((\ga^{-1} f)^{*}*\va)$.
So
$I_{F}(\ga \va)=\ga I_{F}( \va)$  which gives (d).

      Let  $y\in G$.  Since  $F$   is  translation  invariant,
$I_{F}(\va)\subset I_{F}(\Delta_{y}\va)$.
Hence         $sp_{F}(\Delta_{y}\va) $

\noindent $\subset      sp_{F}(\va)$.
Conversely,    suppose   $\ga  \not\in  sp_{F}(\Delta_{z}\va)$
for   all
$z\in  G$,  and  $\ga\not=\hat e$.  Then there exists    $z\in
G$   such  that  $\ga (z)\not=1$   and
$f\in    I_{F}(\Delta_{z}\va)$       such       that     $\hat
{f}(\ga)\not=0$. If   $g=\Delta_{-z} f$  then
$\hat   {g}(\ga)=(\ga  (-z)-1)\hat  {f}(\ga)\not=0$    whereas
$g^{*}*\va=f^{*}*\Delta_{z}\va=0$.   So $g\in I_{F}(\va)$
which implies  $\ga \not\in sp_{F}(\va)$.  This proves (e).

      For  (f),  let $\va$  be constant and  $f\in  F$.   Then
$f^{*} *\va=\va \overline{\hat {f}(\hat e)}$.   So,
if   $\va\not=0$   then $I_{F}(\va)=I(\hat{e},F)$,   which  by
condition (1.5) is dense in  $I(\hat e)$.
Hence   $sp_{F}(\va)=\text  {hull}\,I_{F}(\va)= \text {hull}\,I(\hat {e})=\{\hat
{e}\}$,  the last equality by Wiener's
tauberian   theorem   [20,theorem 7.2.4].    Similarly,     if
$\va=0$    then
$I_{F}(\va)=F$    which   is    dense    in  $L^{1}(G)$.    We
conclude    that
$sp_{F}(\va)=\text {hull}\,L^{1}(G)=\emptyset$.

Finally, $I_{F_{1}}(\va)\subset I_{F}(\va)$   from which
we get (g).
\enddemo

  Following   Reiter[19,p.83]  we  call  a  function     $w\in
L^{\infty}_{loc}(G)$    a weight function on  $G$  if

(2.1)     $w(x)\ge 1$  for all  $x\in G$,  and

(2.2)     $w(x+y) \le w(x)w(y)$  for all  $x,y\in G$.

\noindent Given a weight on  $G$  we define

(2.3)     $L^{1}_{w}(G)   =   \{f  \in   L^{1}(G):  \|f\|_{1,w}=   \int_{G}
|f(x)|w(x)\,dx< \infty\}$.

\noindent Then  $L^{1}_{w}(G)$ is a subalgebra of  $L^{1}(G)$  which is
a Banach algebra under  the norm $\|.\|_{1,w}$.
The Banach space dual of $L^{1}_{w}(G)$   is

(2.4)   $L^{\infty}_{w}(G)=\{\va\in L^{\infty}_{loc}(G):\|\va\|_{\infty,w}=\text{
ess\,sup }_{x\in G} \frac{|\va (x)|}{w(x)}< \infty\}$.

\noindent  As   noted   by   Reiter [19,p.84]   we    may    assume
that   $w$  is upper semicontinuous, and so each  continuous
function  $\va \in L^{\infty}_{w}(G)$  satisfies

(2.5)  $|\va (x)|\le  \|\va\|_{\infty,w} w(x)$  for all  $x\in
G$.

\proclaim {Proposition 2.2} Let  $F=L^{1}_{w}(G)$ or $C_{c}(G)$,  where  $w$
is a weight function on  G.
Then   $F$   is a dense subalgebra of  $L^{1}(G)$   satisfying
(1.3) - (1.5).
\endproclaim

\demo  {Proof}  Since   $C_{c}(G)$  is  dense  in   $L^{1}(G)$
([19,p.68]), so  is   $L^{1}_{w}(G)$.
Property  (1.3) is clear for $F=C_{c}(G)$ and follows from (2.2)
for  $F=L^{1}_{w}(G)$. Property (1.4) is obvious.  For  (1.5),  let
$f\in  I(\ga)$  where  $\ga\in \widehat G$.  Let   $\{f_{n}\}$
be a sequence in  $C_{c}(G)$  converging  to
$f$    in    $L^{1} (G)$.   Take   $u\in C_{c}(G)$  satisfying
$\int _{G} u(t)\,dt=1$,   and   set
$g_{n}  =f_{n}  -  \hat{f} (\ga) \ga  u$.   Then   $g_{n}  \in
I(\ga,C_{c} (G))\subset I(\ga,L^{1}_{w} (G))$  and $\{g_{n}\}$
converges to  $f$  in  $L^{1} (G)$,  proving (1.5).
\enddemo

Let  $F=L^{1}_{w}  (G)$   for  a  weight  function   $w$.   If
$\va\in L^{\infty}_{w} (G)$   and   $f\in F$
then,  for  $y\in G$
$\|f_{y}\va\|_{1}   \le   w(-y)\|f\|_{1,w} \|\va\|_{\infty,w}$
and $\|f\va_{y}\|_{1} \le w(y)\|f\|_{1,w} \|\va\|_{\infty,w}.$
So   $f^{*}*\va$  is well-defined for all  $f\in F$.  In fact,  $f^{*}*\va \in
L^{\infty}_{w} (G)$. We  will
write

(2.6)                $I_{w}(\va)=I_{F}(\va)$               and
$sp_{w}(\va)=sp_{F}(\va).$

 \noindent If in addition  $w$  is bounded, so that  $F=L^{1} (G)$,   we
will write

(2.7)  $I(\va)=I_{F}(\va)$  and  $sp(\va)=sp_{F}(\va).$

\noindent An important additional condition satisfied  by  some weights
is  the
following    condition    of    Beurling-Domar    ([19,p.132],
[5,theorem1.53]):

(2.8)    $\sum^{\infty}_{m=1}\frac{\log   w(mx)}{m^2}<\infty$
for all  $x\in G$.

\proclaim {Proposition 2.3} Let  $w$  be a weight function  on
$G$   satisfying  the
  condition   (2.8)    and    let     $\va\in
L^{\infty}_{w}(G)$.

(a)   $sp_{w}  (f^{*}*\va)  \subset \text{  supp  }\hat  {f}\cap
sp_{w}(\va)$  for all  $f\in L^{1}_{w} (G)$.

(b)   Given a neighbourhood  $V$  of a compact set   $W$    in
$\widehat G$,   there  exists
   $f\in  L^{1}_{w} (G)$   such that $\hat f=1$  on  $W$   and
supp ${\hat {f}} \subset V$.

(c)   If   $f\in L^{1}_{w} (G)$     and  $\hat  f=1$    on   a
neighbourhood  of   $sp_{w} (\va)$   then $f^{*}*\va=\va$.

(d) $sp_{w} (\va)=\emptyset$   if and only if  $\va=0$.
\endproclaim
\demo{Proof}  Since   $L^{1}_{w} (G)$   is  a  Wiener  algebra
([19,p.132]) the proofs  can
be   carried   out  in  the  same  manner   as   for    $L^{1}
(G)$.        See        [6,p.988]      or      [19,p.140-141].
\enddemo
 The following  example   shows that the  last  conclusion  in
Proposition   2.3   is    false  without  the   Beurling-Domar
condition.
\proclaim{Example 2.4} Define a symmetric weight  function   $w$
on  the  discrete
group  of  integers   $\Z$   by   $w(n)=2^{n}  +2^{-n}$.    Let
$\va\in L^{\infty}_{w} (\Z)$   be  given  by
$\va   (n)=2^{n}$      for    $n\in   \Z$.      Then    $sp_{w}(\va)=\emptyset$.    Indeed,
define
$f :\Z\to \r$  by  $f (-1)=2$, $f (1)=-\frac  1{2}$
and $f (n)=0$   if  $n\not=1, -1$.  It
is   easy   to  check  that $f\in L^{1}_{w} (\Z)$ and  $f^{*}*\va=0$, and    hence
$f \in I_{w} (\va)$.    Now
$\widehat {\Z} =T=\r/2\pi \Z$  and for   $t\in T$,
$\hat {f} (t)=2e^{it} - \frac 1{2}e^{-it}$,
which is never  0.  So  $sp_{w}(\va)=\emptyset$.
\endproclaim

\proclaim{Remarks 2.5}

(a)  Let  $w_{1}$ and $w_{2}$    be  two  weight  functions   on
$G$      satisfying    the        condition (2.8).
If       $\va\in  L^{\infty}_{w_{1}} (G)  \cap  L^{\infty}_{w_{2}}  (G)$
then     $sp_{w_{1}}(\va)=sp_{w_{2}}(\va)$.  Indeed,   $w=w_{1}
w_{2}$   is also a weight function on   G
      satisfying the  (2.8).  Since    $w_{1}\le w$   and  $ w_{2} \le
w$
         it     follows     from    [19, p.138]      that
$sp_{w_{1}}(\va)=sp_{w}(\va)=sp_{w_{2}}(\va)$.

(b)   If  $\va\in L^{\infty} (G)$   and   $w$   is  a   weight
function   satisfying
   (2.8), then  $\va\in  L^{\infty}_{w}
(G)$     and so  $sp_{w}(\va)=sp(\va)$.
\endproclaim

Now we
 relate  spectra
relative  to  Beurling  algebras  with  supports  of  Schwartz
distributions.   Let
$S(\r^{n})$   be  the  Schwartz space  of  rapidly  decreasing
infinitely differentiable
complex valued functions on  $\r^{n}$.  Let  $S'(\r^{n})$   be
the dual space of tempered
distributions.  Let  $\va \in L^{\infty}_{loc}(\r^{n})$     be
such that  $T_{\va} (f)=\int_{\r^{n}} \,  f(t)\va(t)\,dt$
   defines  a   distribution   $T_{\va} \in  S'(\r^{n})$.  Then
$\widehat {T}_{\va}(g)=T_{\va} (\hat {g})$   for
$g\in  S(\r^{n})$,  defines the Fourier transform   $\widehat
{T}_{\va}$ of $T_{\va}$.  (See [22, p.146-152]).    As   is
customary, we identify  $\widehat {\r^n}$   with $\r^n$ so that
$\widehat {\r^n}$
becomes an  additive group  with  unit 0.
\proclaim{Proposition 2.6} Let  $w$  be a weight  function  on
$\r^{n}$  satisfying  the
Beurling-Domar   condition  (2.8) and  suppose    $S(\r^{n})\subset
L^{1}_{w}  (\r^{n})$.   If   $\va\in L^{\infty}_{w}  (\r^{n})$
then
$sp_{w}(\va)= \text { supp }{\widehat {T}_{\va}}$.
\endproclaim
\demo{Proof} Let   $\la_{0}\in \r^{n}\setminus sp_{w}(\va)$.
Let  V  be a neighbourhood of   $\la_0$   with
$V\cap  sp_{w}  (\va)=\emptyset$.  By  Proposition 2.3   there
exists    $f\in  I_{w}(\va)$   such  that $\hat  {f}=1$   on a
neighbourhood   $W$  of  $\la_0$    and  $\text { supp }\hat{f}\subset
V$.   Let   $g\in S(\r^{n})$
with  $\text {supp } g\subset  W$.   We  show   $\widehat  {T}_{\va}  (g)=0$.
Indeed     $\overline{\hat  f}g=g$     so
 $ \overline {\hat{g}}^{*}=\overline {\hat{g}}^{*}*f^{*}$.  Hence
$\widehat       {T}_{\va}(g)=T_{\va}      (\hat  {g})=\overline{\hat
{g}}^{*}*\va(0)=\overline {\hat {g}}^{*}*f^{*}*\va(0)=0$,  showing   $\la_0
\not\in \text {supp  } {\widehat {T}_{\va}}$.

     Conversely, suppose $\la_0 \not\in \text {supp  } {\widehat
{T}_{\va}}$.  Let  $V$  be a neighbourhood of  $\la_0$    with
$V\cap \text {supp  } {\widehat {T}_{\va}} =\emptyset$ and
   let $g\in S(\r^{n})$   with  $\text { supp  } g \subset  V$
and   $g(\la_0 )\not=0$.   Let
$f=   \overline   {\hat  {g}}\in  S(\r^{n})\subset   L^{1}_{w}
(\r^{n})$    and   for   a   given    $x\in   \r^{n}$      set
$h(\la)=e^{ix\la} g(\la)$.  Then
$f^{*}*\va(x)=\int_{\r^{n}}\hat {g}(y-x)\va(y)\,dy=T_{\va}      (\hat
{h})=\widehat {T}_{\va} (h)=0$.    Hence    $f\in  I_{w}(\va)$
and  $\hat  {f}(\la_0  )=\overline  {g}(\la_0  )\not=0$.    So
$\la_{0}                 \not\in                 sp_{w}(\va)$.
\enddemo

Next, for a  function   $\va  \in  L^{\infty}
(\r^{n})$   we  define  an  indefinite
integral  $P\va$  by

(2.9)     $P\va   (x_{1} ,  ...,x_{n}   )=\int^{x_{1}}_{0}    \va    (t
,x_{2},...,x_{n})\,dt.$

\noindent It          follows         that          $|P\va        (x_{1}
,...,x_{n})| \le \|\va\|_{\infty} |x_1|$     (a.e.).     Moreover,
$w(x_{1} ,...,x_{n} )=1+|x_1 |$  defines a weight function   $w$
on  $\r^{n}$  satisfying
 (2.8).      Clearly       $\va,P\va\in
L^{\infty}_{w} (\r^{n})$ and
 we prove
\proclaim{Proposition 2.7} If   $\va \in L^{\infty}  (\r^{n})$
and $w(x_{1} ,...,x_{n} )=1+|x_1|$ then $sp_{w}(\va)\subset  sp_{w}(P\va)\subset
sp_{w}(\va)\cup  \{0\}$.
\endproclaim
\demo{Proof}
Consider the tempered distributions $T_{\va}$ and $T_{P\va}$. We have $\frac { \partial
T_{P\va}}{\partial x_{1}} = {T_{\va}}$ and  hence  $\text { supp } \widehat {T}_{\va} \st
\text { supp } \widehat {T}_{P\va}   \st \text { supp } \widehat {T}_{\va}\cup \{0\}$. Applying
Proposition 2.6, we get the required result.
\enddemo

Thirdly,  we discuss spectra  relative  to  the
algebra $F=C_{c}(G)$  and
make  comparisons  with  the spectra  determined  by  Beurling
algebras.  If
   $\va\in  L^{\infty}_{loc}
(G)$ then   $f^{*} *\va$   is  well-defined   and
continuous  for  all   $f \in F$ and   we  will  write

(2.10)  $  I_{c} (\va)=I_{F} (\va)$  and  $sp_{c}  (\va)=sp_{F}(\va)$.
\proclaim{Proposition  2.8}  If    G     is    compact    then
$sp_{c} (\va)=sp(\va)$  for  all
$\va\in \L^{\infty} (G)$.
\endproclaim
\demo{Proof}  By  theorem 2.1(g)  we  have   $sp  (\va)\subset
sp_{c}(\va)$.   Conversely,  if
$\ga  \in \widehat {G}\setminus sp(\va)$  then there exists  $f\in  I(\va)$
with  $\hat {f}(\ga)\not=0$.  Let  $g=f*\ga$  so that
$g\in   C_{c}   (G)$    and   $g^{*}*\va=0$.    Hence    $g\in
I_{c}(\va)$           and          $\hat         {g}(\ga)=\hat
{f}(\ga)\hat{\ga}(\ga)\not=0$.  So
$\ga                   \not\in                   sp_{c}(\va)$.
\enddemo

Example 2.4 above shows also that $sp_{c}(\va)$ may be empty for non-zero $\va$.

 Next,   we   examine    $sp_{c}(\va)$    when
$\va\in L^{\infty}_{loc} (\r)$.
\proclaim{Proposition  2.9} If $\va\in L^{\infty}_{loc}  (\r)$
then  $sp_{c}(\va)= \r$ or   D,   where   D
is  a  discrete  set such that  $\sum \frac{1}{|\la|^{\alf}}$
converges  for  all   $\alf >1$,   the  sum
being over  $\la \in  D \setminus \{0\}$.
\endproclaim
\demo{Proof} If  $I_{c}(\va)=\{0\}$  then   $sp_{c} (\va)=\r$.
On  the  other  hand,  if
$f\in  I_{c}(\va)\setminus \{0\}$  then   $f\in  C_{c}(\r)$.    Hence
$\hat f$   is  an  entire  function  of
exponential  type whose zero set  $Z(\hat f)$   is  such  that
$\sum  \frac{1}{|\la|^{\alf}}$  converges  for  all    $\alf
>1$,  the  sum
being  over  $\la \in Z(\hat f) \setminus \{0\}$. (See Young [23, p.63]).
As   $sp_{c}(\va)\subset  Z(\hat  f)$,   the  result  follows.
\enddemo
Similarly, one can show that
if $\va\in L^{\infty}_{loc}  (\Z)$
then  either   $sp_{c}(\va)=\widehat {\Z} =T$    or
$sp_{c}(\va)$ is finite.
 In contrast,  we  provide the following example  of  a   function
$\va\in L^{\infty} (G)$   with
$sp_{c} (\va)$  neither discrete nor $\widehat {G}$.

\proclaim{Example 2.10} Let $G=\Z\times K$  where   K   is  a  compact
abelian  group.
Define  $\va:G\to \cc$  by $\va (n,k)=\eta (n)\psi (k)$, where
$\eta \in L^{1} (\Z)$  and  $\psi :K\to \cc$
is continuous.  Then  $sp_{c}(\va)=sp_{c} (\eta)\times sp(\psi)$. In
particular, we may choose
$\eta$  and $\psi$  so that  $sp_{c} (\eta)=\widehat{\Z}$   and
$sp(\psi)\not= \widehat K$
\endproclaim

     As an application of  results  in  this  section
we  prove two lemmas needed in section 3.
Firstly,  we  define
another  indefinite integral  $Kf$,  this time for  $f\in
C_{c}(\r)$,  by

(2.11)     $K f(x) = \int_{ -\infty}^{x}  f(t)\,dt.$

\proclaim{Lemma  2.11} Let  $f\in C_{c}(\r)$,     and  let  the
derivatives  $\hat {f}^{(j)}(0)=0$   for
$0\le j\le n-1$  and  $\hat {f}^{(n)}(0)\not=0$ where $n\ge 1$.   Then
$K^{n} f \in C_{c}(\r)$, $\widehat {K^{n}f}(0)\not=0$
and $(K^{n}f)^{(n)}   =f$.
\endproclaim
\demo{Proof} From the assumptions, $\hat {f}(\la)=\la^{n} g(\la)$  for  $\la \in
\r$, where
g  is
entire  and  $g(0)\not= 0$.  Now   $(Kf)'=f$   and  so   $(\widehat {Kf})
(\la)=-i\la^{n-1 } g(\la)$.
Since   $\hat {f}(0)=0$  we get  $Kf\in C_{c} (\r)$.   Arguing  the  same   way,
we  obtain $K^{n} f\in C_{c} (\r)$, $(\widehat {K^{n} f}) (\la)=(-i)^{n} g(\la)$ and
$(K^{n}f)^{(n)}   =f$.
\enddemo
\proclaim{Lemma 2.12} Let $\va \in L^{\infty}_{loc}(\r)$       with      $sp_{c}(\va)$
discrete.    If
$\la \in  sp_{c}(\va)$  then there exists  $h\in C_{c}(\r)$   such that
$h^{*}*\va=\ga_{\la} p,$
where   p is a non-zero polynomial and $\ga_{\la}(t)=e^{i\la t}$.
\endproclaim
\demo{Proof} By theorem 2.1(d),  $0 \in sp_{c}(\ga_{\la}^{-1}\va)$.  Let  $f\in
I_{c}(\ga_{\la}^{-1}\va)$ where
$f\not=0$.
So   0   is  a zero of $\hat f$  of  finite  order   n   say.    Let $g=K^{n} f$.   From
lemma 2.11,  $g\in C_{c} (\r)$, $\hat {g}(0)\not=0 $ and $g^{(n)} =f$.  Hence  $
g\not\in I_{c}(\ga_{\la}^{-1}\va)$
and
therefore $p=g^{*} *\ga_{\la}^{-1}\va\not=0$.  However, $p^{(n)}   =(g^{*}
)^{(n)}*\ga_{\la}^{-1}\va=(-1)^{n} f^{*}*\ga_{\la}^{-1}\va=0$,
 showing that  p  is a polynomial of degree at most  n. Setting
$h=\ga_{\la}g$  we
obtain    $h^{*} *\va=\ga_{\la}g^{*} *\va= \ga_{\la}(g^{*} *\ga_{\la}^{-1}\va)=
\ga_{\la}p$,  as   required.
\enddemo

\head{\bf\S3. Polynomials and their Spectra}\endhead
  In this section we define polynomials on  $G$  and prove the
equivalence
of our definition with that of Domar [5, definition 2.3.2].  We
explore  the
properties  of  polynomials,  and  prove  in  particular   that
$sp_{c} (p)=\{\hat e\}$   for
each non-zero polynomial  $p$.

     A continuous function  $\va:G\to\cc$   will  be  called a
polynomial  of
degree  n  if

(3.1)   $\Delta_{y_{1}}   ...\Delta_{y_{n+1}} \va = 0$  for  all
$y_{1} , ... y_{n+1}    \in G$,  and

(3.2)   $\Delta_{y_{1}}   ...\Delta_{y_{n}} \va\not  =  0$   for
some   $y_{1} , ... y_{n}    \in G$.

  Denote  by  $\Z_{+}$   the semigroup of non-negative integers
with the discrete topology . Set $\Z_{+}^{k}=(\Z_{+})^{k}$, $k\ge 1$.
For arbitrary topological semigroups  $H$   denote by
$C(H$)   the
continuous complex valued functions on  $H$.  The following  two
propositions
will  be  found useful.  The first is proved by induction  and
the  second,  a
simple chain rule for difference operators, is obvious.

\proclaim{Proposition  3.1}  If $  \va\in  C(G)$   then   $\va
(x+my)=\sum^{m}_{j=0} {m \choose j} \Delta^{j}_{y} \va (x)$  for  all
$m\in \Z_{+}$   and  $x,y \in G$.
\endproclaim

\proclaim{Proposition 3.2} Given $\Phi  \in C(G^{n} )$, define
$\Psi \in C(\Z^{n}_{+} \times G^{n})$     by

$\Psi(t_{1}  ,...,t_{n} ,y_{1} ,...,y_{n}  )=\Phi(t_{1}  y_{1}
,...,t_{n}  y_{n} )$, where  $t_{j}\in \Z_{+}$     and   $y_{j}
\in  G$.

 \noindent  If $\Delta_{y,j}$     denotes difference by   $y$
in position  $j$  only, then
$\Delta_{y_{j},j}\Phi(t_{1}  y_{1}   ,...,t_{n}   y_{n}   )=
\Delta_{1,j}\Psi(t_{1} ,...,t_{n} ,y_{1} ,...,y_{n} )$.
\endproclaim

Recall  that  a  function  $\va \in C(G)$  is a polynomial  of
degree  n  in  the
sense  of  Domar  if  $\va(x+ty)$  is a polynomial  in   $t\in
\Z_{+}$   of degree at most  n
for each  $x,y\in G$  and of degree  n  for some  $x,y\in G$.

\proclaim{Proposition 3.3} For  a  function     $\va\in  C(G)$
the   following   are equivalent:

(a)  $\Delta_{y_{1}}   ...\Delta_{y_{n+1}} \va  =  0$   for  all
$y_{1} , ... y_{n+1}    \in G$.

(b) $\Delta_{y}^{n+1} \va = 0$   for all  $y\in G$.

(c)  $\va(x+ty)$  is a polynomial in  $t\in \Z_{+}$   of degree
at  most   n   for  all $ x,y\in G$.

(d)  $\va(x+t_{1} y_{1}+ ...+t_{k} y_{k} )$   is a  polynomial
in  $(t_{1} ,...,t_{k} )\in \Z_{+}^{k}$   of  degree
      at  most  n  for all  $x,y_{1} ,...,y_{k} \in G$,\,  and
$k>0$.
\endproclaim
\demo{Proof} That (a) implies (b) is obvious.  If (b) holds then by
Proposition 3.1  we
have    $\va  (x+ty)=\sum^{n}_{j=0} {t \choose j} \Delta^{j}_{y}  \va  (x)$,
showing   that   (b) implies (c).     Given (c)    and  $x,y_{1}
,...,y_{k}  \in  G$,  define     $\psi(t)=\va(x+t_{1}   y_{1}+
...+t_{k}  y_{k} )$    for $(t_{1} ,...,t_{k} )\in \Z^{k}_{+}$.
Then   $\psi$  is a polynomial in each  $t_{j}$    of   degree
at  most  n.  So $\psi$   is a polynomial in   t   of   degree
at  most   kn.   Hence

$\psi(t)=  \sum_{|\alf|\le kn}   a_{\alf} t^{\alf}$,  where
$a_{\alf}\in     \cc$,    $\alf=(\alf_{1},...,\alf_{k})\in
\Z^{k}_{+}$ and   $|\alf|= \alf_{1}+...+\alf_{k}$.

\noindent Moreover,   for    $\la \in \Z_{+}$,  the function $\psi(\la t)= \sum_{|\alf|\le
kn}    a_{\alf}\la^{|\alf|} t^{\alf}=\va(x+\la t_{1} y_{1}+
...+\la t_{k} y_{k})$    is
again  a  polynomial in   $\la$   of  degree   at   most    n.
Hence    $a_{\alf}=0 $  for
$| \alf| > n$,   showing that (c) implies (d).
 Finally,   given
(d)  and  $x,y_{1} ,...,y_{n+1} \in G$,
define $\psi(t)=\va(x+t_{1} y_{1}+ ...+t_{n+1} y_{n+1})$   for
$(t_{1} ,...,t_{n+1} )\in \Z^{n+1}_{+}$.   Then
$\psi$   is a polynomial in  t  of degree at most  n.  So  all
(n+1)th  differences
of  $\psi$  are zero.  Hence, by Proposition 3.2,
  $\Delta_{y_{1}}    ...\Delta_{y_{n+1}}  \va(x)=\Delta_{1,1}
...\Delta_{1,n+1} \psi(0)=0$,  showing that (d) implies  (a).
\enddemo

\proclaim{Theorem 3.4} For  a  function    $\va\in C(G)$   the
following   are equivalent:

(a)  $\va$  is a polynomial of degree  n.

(b)  $\Delta_{y}^{n+1}  \va  =  0$    for  all   $y\in  G$  and
$\Delta_{y}^{n} \va \not= 0$ for some  $y\in G$.

(c)  $  \va$  is a polynomial  of degree   n in the  sense  of
Domar.

(d)  $\va(x+t_{1} y_{1}+ ...+t_{k} y_{k} )$   is a  polynomial
in  $(t_{1} ,...,t_{k} )\in \Z^{k}_{+}$   of  degree
      at  most  n  for all  $x,y_{1} ,...,y_{k} \in G,  \,k>0$
and of degree n for some  $x,y_{1} ,...,y_{k} \in G, \,k>0$.
\endproclaim
\demo{Proof}  Given (d), we have Proposition  3.3  (c).
Moreover we  can
choose  $x,y_{1} ,...,y_{k}$    such that  $\va(x+t_{1} y_{1}+
...+t_{k}  y_{k}  )$   is a polynomial in  $(t_{1}  ,...,t_{k}
)\in  \Z^{k}_{+}$    of   degree n.   So   $\va(x+t_{1}  y_{1}+
...+t_{k}    y_{k}   )=   \sum_{|\alf|\le   n}     a_{\alf}
t^{\alf}$,  where     $a_{\alf}\in \cc$, and $\sum_{|\alf|=
n}  a_{\alf} t^{\alf}$  is not identically zero.  Choose   t
so  that  $\sum_{|\alf|=  n} a_{\alf}  t^{\alf}\not=0$  and
set      $y=t_{1} y_{1}+ ...+t_{k} y_{k}$.        Then
$\va (x+\la y)=\sum_{|\alf|\le n}   a_{\alf} t^{\alf}\,
\la^{|\alf|}$      is  a  polynomial  in  $ \la \in \Z_{+}$    of
degree   n, showing that (d) implies (c).
      Given(c), we have Proposition 3.3(b).   Moreover,
there  exist $x,y\in G$   such  that   $\va (x+ty)$   is    of
degree     n     in    $t  \in  \Z_{+}$.     So  $\va   (x+ty)=
\sum^{n}_{j=0} a_{j} t^{j}$
   where      $a_{j}   \in   \cc$,   $a_{n}\not=0   $.       By
Proposition      3.2,     $\Delta_{y}^{n}\va(x)=\sum^{n}_{j=0}
a_{j}\Delta_{1}^{n} \, t^{j}=a_{n}  n  !\not=0$,   showing that   (c) implies (b).
  That  (b)  implies (a)   is obvious.   Finally,  if  (a) holds  then   Proposition  3.3(d)
holds.   Moreover,
there   exist   $x,y_{1}  ,...,y_{n}\in  G$       such    that
$\Delta_{y_{1}}   ...\Delta_{y_{n}} \va(x)\not=0$. So
  $\va(x+t_{1}  y_{1}+ ...+t_{n} y_{n} )= \sum_{|\alf|\le  n}
a_{\alf} t^{\alf}$  and       by       Proposition 3.2
$\Delta_{y_{1}}    ...\Delta_{y_{n}} \va(x)=a_{\beta}$     where
$\beta=(1,1,...,1)\in  \Z_{+}^{n}$.   So    $a_{\beta}\not=0$
 and $\va(x+t_{1} y_{1}+ ...+t_{n} y_{n})$
 is a  polynomial  in   $(t_{1} ,...,t_{n} )$   of  degree n,
showing          that         (a)         implies         (d).
\enddemo

\proclaim{Remarks 3.5}

(a)   On any abelian group  G  the polynomials of degree 0  are
precisely the     constant functions.

(b)    Polynomials of degree  larger  than  0  are   unbounded.
So   the   only      polynomials on a compact  group  are  the
constant functions.

(c)   The definition of polynomials and the results proved  so
far  in   this   section require only that  G  be  an  abelian
topological semigroup.

(d)   Any  polynomial   $\va$   on    G    of   degree   1   is
of   the   form
        $\va     (x)=\va    (0)+\alf (x)$,     where     $\alf (x+y)=\alf (x)+\alf (y)$  for
$x,y\in G$.   So   $\va$
is uniformly continuous on  $G$.

(e)    From  theorem 3.4 we see that any polynomial  $\va$   of
degree N   satisfies

  $\lim_{|n|\to \infty} \frac{|\va(nx)|}{|n|^{N+\alf}}$
 $=0$ for all  $x\in G$  and all  $\alf >0$.

 \noindent  Conditions similar
     to this will be needed in section 4.
\endproclaim

We turn now to the investigation of spectra of polynomials.

\proclaim{Proposition 3.6} Let  p  be a  non-zero   polynomial
on   G.   Suppose
$f^{*}*p$  is well-defined for all  $f\in F$, where  $F$  is a
dense  subalgebra  of
$L^{1}  (G)$   satisfying (1.3) - (1.5).  Then   $sp_{F}  (\ga
p)=\{\ga\}$  for all  $\ga \in \widehat G$.
\endproclaim
\demo{Proof}  Let   p   have  degree  n.   By   theorem   3.4,
$\Delta_{y_0}^{n}p\not=0$   for  some
$y_{0} \in G$.     Applying     theorem 2.1(e)      n     times
we     obtain
$sp_{F}   (\Delta_{y_0}^{n}p)\subset   sp_{F}(p) \subset    \{\hat
e\}\cup_{y\in  G}  sp_{F}  (\Delta_{y}^{n}  p)$.     But    by
theorem   2.1    (f), $sp_{F} (\Delta_{y_0}^{n}p) =\{\hat e\}$
and  $sp_{F} (\Delta_{y}^{n}p) \subset \{\hat e\}$    for each
$y\in G$.  Hence  $sp_{F} (p)= \{\hat e\}$
and    by    theorem   2.1(d),    $sp_{F}   (\ga    p)=\{\ga\}$.
\enddemo

Immediately we obtain

\proclaim{Corollary 3.7} Let  p  be a non-zero polynomial   on
G   and  suppose $\ga\in  \widehat G$. Then

(a)  $sp_{c} (\ga p)=\{\ga\}$,

(b)  $sp_{w} (\ga p)=\{\ga\}$ for any weight  $w$  such that $p\in
L_{w}^{\infty} (G)$.
\endproclaim
As an application of the above corollary we show
\proclaim {Theorem 3.8} Let  $w$  be a weight function on  $\r$  and   let $\va\in
L^{\infty}_{loc}  (\r)$.
If $sp_{c}(\va)$  is dicrete then $sp_{w}(\va)=sp_{c}(\va)$ .
\endproclaim
\demo{Proof} By   theorem 2.1(g)   we  have  $sp_{w}(\va)\st sp_{c}(\va)$.
Conversely,  let $\la \in sp_{c}(\va)$.
 If $\la\not \in sp_{w}(\va)$      then  there  exists    $g \in I_{w}(\va)$ such  that $\hat
{g}(\la)\not=0$.   By    lemma 2.12    there    exists $h\in C_{c} (\r)$
such   that
$h^{*}*\va=\ga_{\la}p$,  where  p  is a non-zero polynomial.  Now   $g^{*} *h^{*}
*\va=0$ so
$g\in I_{w}  (h^{*}  *\va)$,   but  by  corollary 3.7  $sp_{w}  (h^{*}  *\va)=\{\la\}$.
Hence
$\hat {g}(\la)=0$,
a contradiction showing  $\la \in sp_{w}(\va)$.
\enddemo
      The  following example shows that for the conclusion  of
theorem 3.8  it
is not enough to know that  $sp_{w}(\va)$  is discrete.

\proclaim {Example 3.9} Consider   the    function     $\va$     on     $\r$
given    by
$\va (t)=\sum_{n=1}^{\infty}n^{-2} e^{i \,t n^{1/2}}$. Then   $\va$   is  almost
periodic  and   it
is  known
([15] or [12,V.5])  that  $sp(\va)=\{ n^{1/2}:n\in N\}$.  We  prove  that    $sp_{c}(\va)=
\r$.
Indeed,   $sp (\va)\st sp_{c}(\va)$    and    if   $sp_{c}(\va) \not=\r$    then
by
Proposition 2.9, $\sum_{n=1}^{\infty} n^{-\alf / 2} < \infty$
  for all  $\alf >1$. As this is false,  $sp_{c}(\va) =\r$.
\endproclaim

\head{\bf\S4. Primary ideals and functions with finite  spectrum}\endhead

     In this section we study Beurling  algebras   $L^{1}_{w} (G)$  when the  weight
function   $w$  satisfies the following two conditions for  some $N\in \Z_{+}$:

(4.1) $\text {sup}_{m\in \Z} \frac {w(mx)}{1+|m|^{N+\alf (x)}} < \infty$  for all  $x\in G$
and some  $\alf (x)<1$;

(4.2) $\text {inf}_{m\in \Z} \frac {w(mx)}{|m|^{N}} > 0 $  for some  $x\in G$.

It   should   be  noted  that  condition  (4.1)  implies   the
Beurling-Domar
condition (2.8).   Using   conditions    (4.1) and  (4.2)    we
characterize  those
functions  $\va\in L^{\infty}_{w} (G)$  for which  $sp_{w}(\va)$  is finite.

      We  begin by studying the primary ideals of $L^{1}_{w} (\Z)$   when
$w$   satisfies
(4.1) and  (4.2) with  $G=\Z$.  Recall that  a  primary   ideal   is
one  that  is
contained in exactly one maximal ideal, and the cospectrum  of
an  ideal   J
of  $L^{1}_{w} (G)$ is the set
$\{\ga \in \widehat {G}: \hat {f}(\ga) = 0\text { for all  } f \in J \}.$

     When  $G=\Z$,  conditions (4.1) and (4.2) are equivalent to the existence
of $ c_{1} >0, c_{2} >0$  and $\alf\in  (0,1)$  such that

(4.3)  $c_{1} (1+|m|^{N} )\le w(m)\le c_{2} (1+|m|^{N+\alf})$       for all  $m\in \Z$.

\noindent {In} this case,   $\{ \hat f:f\in L_{w}^{1} (\Z) \}$   is   a subalgebra of  $C^{N} (T),$  the
algebra
of   continuous  functions  on   $T=\r / 2\pi \Z$  whose derivatives up
to order  N
are also continuous on  $T$.  We then define

(4.4) $I_{k} =\{ f\in L_{w}^{1} (\Z):\hat {f}^{(j)}(\hat e) = 0    \text { for } 0 \le j \le
k\}$, where  $0\le k\le N$.

\noindent {So}   $I_{k}$     is   a  closed  primary  ideal  of $L_{w}^{1} (\Z)$   with
cospectrum  $\{\hat{e}\}$.

\proclaim{Theorem 4.1} Let   $w$   be  a  weight  function   on     $\Z$
satisfying  (4.3).  Then   $I_0 ,I_1 ,...,I_N$    are  the  only  closed
primary
ideals of  $L_{w}^{1} (\Z)$   with cospectrum  $\{\hat{e}\}$.
\endproclaim
\demo{Proof}  Let  $e_{m} (n)=\delta_{m,n} $ denote     the  Kronecker  function,
where
$m,n \in \Z$.
Then  $e_{1}$  and $e_{-1}$    generate the algebra  $L_{w}^{1} (\Z)$,  and   $\hat
{e_{1}} (t)=e^{it}$  and  $\hat {e_{-1}} (t)=e^{-it}$
for   $t\in T$.   Moreover,  $e_{1}^{m} =e_{m}$   so  $\alf_{m}=
\|e_{1}^{m}\|_{1,w}=w(m)$   for   each
$m \in \Z$. It follows from (4.3)   that
 $\lim_{r \to 1-}
(1-r)^{N+2}\sum_{m=0}^{\infty}\alf_{m}r^{m}=\lim_{r \to 1-}
(1-r)^{N+2}\sum_{m=0}^{\infty}\alf_{-m}r^{-m}=0$.
 By a theorem of  Gelfand  [8]  or   [16,
p.209], $L_{w}^{1} (\Z)$ has   exactly    N+1    ideals    with    cospectrum
$\{\hat{e}\}$,    namely  $I_0 ,I_1 ,...,I_N.$                 \enddemo

\proclaim{Proposition 4.2} Let   $w$   be  a  weight  function  on   $\Z$
satisfying (4.3).  If  $0\le k \le N$  then  $\Delta_{m}^{k+1} f\in I_{k}$   for   all
$m\in \Z$
and all  $f \in L_{w}^{1} (\Z)$.
\endproclaim
\demo{Proof} If  $ g= \Delta_{m}^{k+1} f$  then $\hat {g}(t)=(e^{itm}   -1)^{k+1} \hat
{f}(t)$.    Hence
$g\in I_{k}$ as   claimed.
\enddemo
\proclaim{Corollary 4.3} Let   $w$   be  a  weight  function  on  $\Z$ satisfying  (4.3).   If
$\va \in L_{w}^{1} (\Z)$   and   $sp_{w} (\va)=\{\hat{e}\}$,   then $\va$
is  a polynomial of degree at most  N.
\endproclaim
\demo{Proof}  Since    $I_{w} (\va)$    is   a   closed   ideal   with
cospectrum  $\{\hat{e}\}$,  we
conclude  from  theorem 4.1 that  $I_{w} (\va)=I_{k}$   for some   k, $0\le k \le N$.
Now for
arbitrary    $f \in L_{w}^{1} (\Z)$    and    $m\in \Z,$     we    have     $\Delta_{-
m}^{k+1}f \in I_{k}$
and   so $f^{*}*\Delta_{m}^{k+1}\va
=(\Delta_{-m}^{k+1}f)^{*}*\va=0$.  So $I_{w}(\Delta_{m}^{k+1}\va)=L_{w}^{1} (\Z)$
and  $sp_{w}(\Delta_{m}^{k+1}\va)=\emptyset$.
By  Proposition 2.3 (d), $\Delta_{m}^{k+1}\va=0$,  and  by  (3.1)   $\va$    is   a
polynomial  of
degree                at                most                k.
\enddemo

     We now consider functions on more general groups  G.

\proclaim{Theorem 4.4} Let  $w$  be a weight function satisfying (4.1) and (4.2)  and
let  $\va \in L_{w}^{\infty} (G)$.  Let  $sp_{w} (\va)=\{ \ga_{1} ,...,\ga_{n} \}$ be  a
finite subset of $\widehat G$.
Then
there  exist polynomials  $p_{1} ,...,p_{n}$   in $ L_{w}^{\infty} (G)$    of  degree  at
most  N   such
that  $\va= \sum_{k=1}^{n}\ga_{k}p_{k}$.
\endproclaim
\demo{Proof}  As $sp_{w} (\va)$  is compact,  $\va$  coincides a.e.  with a continuous
function.  (See [19, p.142]).  We  may  therefore  assume
$\va\in C(G)$.
Now  suppose  $sp_{w} (\va)=\{ \hat{e}\}$.  For  $y\in G$,  let  $H_{y}$   be the
subgroup
generated
by   y   and  let  $K_{y}$   be  its closure.   Let   $w_{y}$    be  the
restriction of   $w$   to
$K_{y}$    and  $\va_{y} (h)=\va (h+y)$  for  $h \in K_{y}$.  Then $\va_{y}\in
L_{w_{y}}^{\infty} (K_{y})$   and by a
theorem
of  Reiter (see [18] or [19, p.144])   $sp_{w_{y}} (\va_{y}) \st \pi (sp_{w} (\va))= \{\hat
{e_{y} }\}$,
where  $\pi : \widehat {G} \to\widehat {K_{y}}$   is the canonical projection and  $\hat
e_{y}$    is  the
unit in  $\widehat {K_{y}}$.
By  a  theorem  of Weil (see [21, p.96-97]  or  [19, p.88])
either  $K_{y}$     is
compact  or   $H_{y}$    is isomorphic to  $\Z$  $(H_{y} =\Z)$.   If   K    is
compact then   $\va$
is  bounded with   $sp(\va_{y})\st \{ \hat e_{y}\}$.  Therefore   $\va_{y}$  is
a  constant.
(See  Rudin  [20,  6.8.3].)   On  the  other  hand, by
corollary 4.3,  if $H_{y} =\Z$
then      $\va$    is a polynomial of degree  at   most    N.
Applying  these
results  to   $\psi,$  defined by  $\psi(z)=\va (x+z)$   for    $x,z\in G,$    and
noting  that
$sp_{w}  (\psi)=sp_{w} (\va)$  we conclude that  $\va (x+my)$  is a polynomial  in
$m\in \Z_{+}$    of
degree at most  N  for each  $x,y\in G$.  So  $\va$  is a polynomial of
degree  at
most  N.
  For  the  general  case,  $sp_{w} (\va)=\{ \ga_{1} ,...,\ga_{n} \},$   choose  a
neighbourhood  V
of  $\{\hat {e}\}$  in  $\widehat {G}$  and functions  $f_{1} ,...,f_{n}\in L_{w}^{1}
(G)$  such  that  $\hat {f}_{k}=1$
on
$\ga_{k}  V$   and    $\ga_{j}  V\cap \ga_{k} V=\text {  supp    } \hat{f_{j}}\cap\text{
supp  }\hat{f_k} =\emptyset$   for   $j,k=1,...,n$
and
$j\not=k$.   (See   Proposition  2.3).   Then    $sp_{w}  (f_{k}^{*}  *\va) \st  \text {
supp }\hat{f_{k}}
\cap sp_{w}(\va)=\{\ga_{k}\}$
and  by  the  first part of this  proof,   $\ga^{-1}_{k}(f_{k}^{*}  *\va)=p_{k}$,   a
polynomial  of
degree  at most  N.  So  $f^{*} *\va= \sum_{k=1}^{n}\ga_{k}p_{k}$,   where
$f=\sum_{k=1}^{n}f_{k}.$  But
 $\hat f=1$
on   a   neighbourhood   of   $sp_{w} (\va)$    and   so,   again   by
Proposition   2.3, $f^{*}*\va=\va$.
\enddemo
\head{\bf \S5. Uniformly continuous functions with bounded differences}\endhead
 Let   $C_{u}  (G)$   be  the  space  of  uniformly  continuous
complex  valued
functions on  G.  Let  $C_{ub}  (G)=C_{u} (G)\cap L^{\infty} (G)$  and denote  by
$C_{ud}(G)$   the
linear  subspace of  $C_{u} (G)$  consisting of functions   $\va\in C_{u}  (G)$
satisfying

(5.1)  $\Delta_{h}  \va \in C_{ub}  (G)$  for all  $h\in G$.

We  notice that  if   $G$   is  connected  then  by  a   theorem
of  Pontrjagin
(see[17, p.377]  or  [19,  p.97])    $G=K\times\r^{n}$,     where    K    is
compact  and
connected and  $n\ge 0$.  It follows that  $C_{ud}  (G)=C_{u} (G)$.

      In  this section we study the spectra of functions    $\va \in C_{ud}(G)$.   To  do
this  we  construct  associated  Beurling  algebras.   Indeed,
80adefine   $w_{\va} :G \to \cc$ by

(5.2) $w_{\va} (y) = 1 + 1/2 (\|\Delta_{y} \va\|_{\infty} + \|\Delta_{-y} \va\|_{\infty})$,
$y \in  G$.

\proclaim{Proposition 5.1} If  $\va \in C_{ud}(G)$  then $w_{\va}$    is  a   uniformly
continuous
weight function on  G  satisfying the condition

(5.3)  $\lim_{|m|\to \infty} \frac{w_{\va}(mx)}{m^{2}}=0$  for all  $x \in G$.

\noindent Moreover  $\va \in L_{w_{\va}}^{\infty}  (G).$
\endproclaim
\demo{Proof}    Let     $\alf (y)= w_{\va}(y)-1$.
 Then    $\alf (y_{1} +y_{2})\le
\alf (y_{1}
)+\alf (y_{2} )$,
$\alf (y)=\alf (-y )$,   and  $|\alf (y_{1} +y_{2})\ - \alf (y_{2})|\le \alf (y_{1} )$    for  $ y,
y_{1}, y_{2} \in G$.
It
follows  that  $\alf$ and $w_{\va}$   are uniformly continuous on  G,   and   that
$w_{\va}$
is  a  weight
function.       Moreover,     for      $n\in \Z$ and        $y \in  G$, $w_{\va}(my)= 1+\alf
(my) \le 1+ |m| \alf (y),$
from        which        (5.3)       follows.         Finally,
$|\va(x)| \le |\va (x)- \va (0)| +|\va (0)|\le 2 \alf (x)+|\va (0)|\le (2+ |\va (0)|)w_{\va}(x)$,
for  $x\in G$,  showing    $\va \in L_{w_{\va}}^{\infty}  (G)$.
\enddemo

      We remark that condition (5.3) implies  condition  (4.1)
with   N=1.
This in turn implies the Beurling-Domar condition (2.8).
    In  the following  we extend  a result  proved by Basit[1] for the case $G=\r$.

\proclaim{Theorem 5.2} If $\va \in C_{ud}(G)$  and  $\hat {e} \not\in sp_{w_{\va}}
(\va)$   then   $\va$  is
bounded.
\endproclaim
\demo{Proof}  Let   V   be  a  compact  neighbourhood   of     $\hat {e}$
such   that
$V\cap sp_{w_{\va}} (\va)=\emptyset$.   By  Proposition 2.3 there exists    $f\in
L_{w_{\va}}^{1} (G)$
such  that
supp $\hat {f}\st V$   and  $\hat {f}(\hat {e})=1$.  So  $sp_{w_{\va}}  (f^{*}
*\va)=\emptyset$    and   hence   $f^{*} *\va =0$.
Hence           $\va (t)=\va (t)-f^{*} *\va (t)=\int_{G}   [\va(t)-\va(t-s)]f^{*}   (s)\,ds$
and $|\va (t)|\le 2\int_{G}   w_{\va}(s) |f^{*}   (s)|\,ds=2\|f\|_{1,w_{\va}}$,  for   $t\in G$.
So   $\va$   is
bounded.
\enddemo
We  return  now  to  the consideration   of   indefinite
integrals   $P\va$   of
functions  $\va \in L^{\infty} (\r^{n} )$.
\proclaim{Theorem 5.3} If  $\va \in L^{\infty} (\r^{n})$, $0 \not\in sp (\va)$
and  $P\va \in C_{u} (\r^{n})$    then
$P\va$
is bounded.
\endproclaim

\demo{Proof}  If   $w(x_{1} ,...,x_{n} )=1+|x_{1}|$    then   $w$   is  a  weight  on
$\r^{n}$
satisfying (2.8)  and   $\va, P\va \in L^{\infty}_{w} (\r^{n} )$.
By  remark
2.4 (b),   $sp(\va)=sp_{w} (\va)$.
Let   V   be a neighbourhood of  $0$    with  $sp(\va)\cap V=\emptyset$.    By
Proposition
2.3 (b),  there  is  a  function   $h \in L_{w}^{1} (\r^{n})$   such  that   $\hat {h} =1$
on   a
neighbourhood     W     of  $0$        and     supp $\hat {h} \st V$.      By
Proposition 2.3(a),
$sp_{w}(h^{*} * P\va)\st \{0\}$    and   $h^{*} *\va=0$.  This implies that $h^**P\va$ is infinitely differentiable. By assumption $P\va$ is uniformly continuous, hence   $h^{*} *P\va=0$ is also uniformly continuous. Since $\frac{\partial}{\partial{x_1}}h^{*} *P\va=0$, $ h^{*} *P\va(x_1,x_2,\cdots,x_n)=\psi(x_2,\cdots,x_n)$. As $P\va \in L_w^{\infty} (\r^n)$, $h^** P\va \in L_w^{\infty} (\r^n)$. This implies that $\sup_{x\in \r^n} |\psi(x_2,\cdots,x_n)| \le  \sup_{x\in \r^n} \frac{|\psi(x_2,\cdots,x_n)|}{1+|x_1|}< \infty$ showing that $\psi$ is bounded. Set $\eta =\va-\psi$. Let $g\in L_w^1 (\r^n)$ such that supp$\hat{g} \st W$ and $\hat{g} (0)\not = 0$. We have $g^**\eta =0$. Hence
$0\not \in sp_w(\eta)$ so by Theorem 5.2, $\eta$ is bounded. Hence $P\va=\eta+\psi$ is bounded.
\enddemo

\proclaim {Corollary 5.4} If  $\va \in L^{\infty} (\r )$ and $0 \not\in sp (\va)$       then
$P\va$
is bounded.
\endproclaim
\demo{Proof} Since $| P\va(x)-P\va(y) | \le ||\va||_{\infty}|x-y|$, $P\va \in C_{u}(\r)$   and
the
corollary        follows       from        the        theorem.
\enddemo
\proclaim{Remarks 5.5}

(a) Let $\va(x)=x e^{ix}$ and $w(x)=1+|x|$ for $x\in \r$. By corollary 3.7
$sp_{w}(\va) =\{1\}.$ So $0\not\in sp_{w}(\va)$ yet $\va$ is unbounded. Hence
 a condition like $\va \in C_{ud}(G)$  is necessary in theorem 5.2.

(b) Let $\va(x,y)= e^{i y}$ for  $(x,y) \in \r^{2}$. Then  $sp(\va)=\{ (0,1)\}.$
So $(0,0)\not\in sp(\va)$  yet $P\va (x,y)=x e^{i y}$ is
unbounded. Hence a condition like $P\va \in C_{u}(\r)$ is necessary in theorem
5.3.

(c)   Favard [7]  obtained corollary 5.4 for almost  periodic
functions   $\va$.
Later, Levitan-Zhikov [13,p.89] obtained it for  $\va \in C_{ub}(\r)$.

(d) We thank D. Domar who kindly informed us that the proof of Theorem 5.3 appearing in Bull. Australian Math. Soc., is not complete. Also he sent us the proof of the following sharper result:

If $\psi$ is uniformly continuous on $\r\times \r^{n-1}$, bounded on $0\times \r^{n-1}$ and if supp $\widehat{\frac{\partial {\psi}}{\partial{x_1}}}= sp (\frac{\partial {\psi}}{\partial{x_1}})$ does not contain $0$, then $\psi$ is bounded (Distribution theory notation).
\endproclaim

Address: School of Mathematical Sci., P.O.Box No. 28M,  Monash University, Vic. 3800, Australia.

 e-mail bolis.basit\@monash.edu,\,  alan.pryde\@monash.edu

\enddocument